\documentclass[12pt]{article}
\usepackage[cp866]{inputenc} % if typed in DOS
\usepackage[english]{babel}
\usepackage{graphicx}
\usepackage{amsmath}
\usepackage{float}

\title{A generalized \((G'/G)\) - expansion method for the loaded modified Korteweg-de Vries equation}

\author{I. I. Baltaeva$^1$, I. D. Rakhimov$^2$, M. M. Khasanov$^3$}
\begin{document}

\maketitle

\emph{$^1$$^,$$^2$$^,$$^3$ Urgench state university, H. Alimdjan 14, Urgench 220100, Uzbekistan}
\begin{center}
\textit{e-mail: iroda-b@mail.ru$^1$, \, ilxom.raximov.87@gmail.com$^2$, \,\\ hmuzaffar@mail.ru$^3$}	
\end{center}

\begin{abstract}
This paper is dedicated to finding the solutions of the equation of the loaded modified Korteweg-de Vries. By the way, it is shown to find the solutions via (G'/G) - expansion method that is one of the most effective ways of finding solutions.
\end{abstract}
\section{Introduction}
It is known that loaded differential equations have great practical applications. In the literature \cite{1,2,3,4,5}, loaded differential equations are typically called equations containing in the coefficients or in the right-hand side any functionals of the solution, in particular the values of the solution or its derivatives on manifolds of lower dimension. These types of equations were explored in the works of N.N. Nazarov and N.N. Kochin. However, they did not use the term ``loaded equation''. At first, the term has been used in the works of A.M. Nakhushev, where the most general definition of a loaded equation is given and various loaded equations are classified in details, for instance, loaded differential, integral, integro-differential, functional equations etc., and numerous applications are described. 

In the past several decades, finding solutions to nonlinear evolution equations has been studied by many researchers. There are direct and inverse methods for finding solutions to integrable nonlinear evolution equations. In particular, the solutions of the integrable nonlinear evolution equations were found by using the Hirota direct method \cite{6,6_1}, inverse scattering problem \cite{7,8,8_1} and Darboux transformation \cite{9,10}. Alternatively, the $(G'/G)$ - expansion method \cite{11,12,18,14,15,16,17,18_1,19,20,20_1} is also effective in finding traveling wave solutions of nonlinear evolution equations. 

Integration of the loaded modified Korteweg -- de Vries (mKdV) equation in the class of periodic functions is studied in \cite{21}. 

In this article, the solutions of the loaded mKdV equation, are explored by usage of $(G'/G)$ - expansion method.

Consider the following loaded mKdV equation
\begin{equation}
\label{1}
q_{t} -6q^{2} q_{x} +q_{xxx} -\gamma (t)q(0,t)q_{x} =0,
\end{equation}
where $q(x,t)$ is an unknown function, $x\in R$, $t\ge 0$, $\gamma (t)$ - is the given real continuous function. 
\section{Description of the generalized $(G'/G)$ -expansion method}

Let us be given a nonlinear partial differential equation in the form below
\begin{equation}
\label{2}
F(q,q_{t} ,q_{x} ,q_{tt} ,q_{xx} ,q_{xt} ,...)=0,
\end{equation} 
with two independent variables \textbf{$x$ }and \textbf{$t$}. \textbf{$q=q(x,t)$ }is a unknown function, \textbf{$F$ }is a polynomial in \textbf{$q$} and its partial derivatives in which the highest order derivatives and nonlinear terms are involved. Now we give the main steps of the $(G'/G)$ -expansion method \cite{22}:

\textit{Step 1.} We use the travelling wave transformation in the following form
\begin{equation}
\label{3}
q(x,t)=q(\xi), \, \xi =kx+\Omega (t),
\end{equation}
where $k$ is a parameter and $\Omega (t)$ is a continuous function which depends on $t$. We reduce equation (\ref{2}) to the following nonlinear ordinary differential equation:
\begin{equation}
\label{4}
P(q,q',q'',q''',...)=0,
\end{equation}
where $P$ is a polynomial of $q(\xi )$ and all its derivatives $q'=dq(\xi )/d\xi $, $q''=d^{2} q(\xi )/d\xi ^{2}$, \dots  .

\textit{Step 2.} We assume that the solution of equation (\ref{4}) has the form:
\begin{equation}
\label{5}
q(\xi )=\sum _{j=0}^{m}a_{j} \left(\frac{G'}{G} \right)^{j}.
\end{equation}
Here $G=G(\xi )$ satisfies the following second order ordinary differential equation
\begin{equation}
\label{6}
G''+\lambda G'+\mu G=0,
\end{equation}
where $G'=dG(\xi )/d\xi $, $G''=d^{2} G(\xi )/d\xi ^{2}$ and $\lambda $, $\mu $, $a_{j}$ $(j=1,2,\, ...\, ,\, m)$ are constants that can be determined later, provided $a_{m} \ne 0$.

\textit{Step 3.} We determine the integer number $m$ by balancing the nonlinear terms of the highest order and the partial product of the highest order of (\ref{4}).

\textit{Step 4.} Substitute (\ref{5}) along with (\ref{6}) into (\ref{4}) and collect all terms with the same order of $\left(\frac{G'(\xi )}{G(\xi )} \right)$, the left-hand side of (\ref{4}) is converted into a polynomial in $\left(\frac{G'(\xi )}{G(\xi )} \right)$. Then, set each coefficient of this polynomial to zero to derive a set of over-determined partial differential equations for $a_{j} $ $(j=1,2,\, ...\, ,m)$ and $\xi $.

\textit{Step 5.} Substituting the values $a_{j} $ $(j=1,2,\, ...\, ,m)$ and $\xi $ as well as the solutions of equation (\ref{6}) into (\ref{5}) we have the exact solutions of equation (\ref{2}).
\section{Exact Solutions of the loaded mKdV equation}

In this section, we will show how to find the exact solution of the loaded mKdV equation using the $(G'/G)$ -expansion method. For doing this, we perform the steps above for equation (\ref{1}). The travelling wave variable below

\begin{equation} 
\label{7}
q(x,t)=q(\xi ), \, \xi =kx+\Omega (t),     
\end{equation} 
permits us converting equation (\ref{1}) into an ordinary differential equation for $q=q(\xi )$

\begin{equation}
\label{8}
\Omega '_{t} (t)q'-6kq^{2} q'+k^{3} q'''-k\gamma (t)q(0,t)q'=0,
\end{equation} 
integrating it with respect to $\xi $ once yields to

\begin{equation}
\label{9}
C+\Omega '_{t} (t)q-2kq^{3} +k^{3} q''-k\gamma (t)q(0,t)q=0,
\end{equation} 
where $C$ is an integration constant that can be determined later.

We express the solution of equation (\ref{9}) in the form of a polynomial in $(G'/G)$ below

\begin{equation}
\label{10}
u(\xi )=\sum _{j=0}^{m}a_{j} \left(\frac{G'}{G} \right)^{j},
\end{equation}
where $G=G(\xi )$ satisfies the second order ordinary differential equation in the form
\begin{equation}
\label{11}
G''+\lambda G'+\mu G=0.
\end{equation} 
Using (\ref{10}) and (\ref{11}), $q^{3}$ and $q''$ are easily derived to

\begin{equation}
\label{12}
q^{3} (\xi )=a_{m}^{3} \left(\frac{G'}{G} \right)^{3m} +\, ... , 
\end{equation} 

\begin{equation}
\label{13}
q''(\xi )=m(m+1)a_{m} \left(\frac{G'}{G} \right)^{m+2} +\, ... .
\end{equation} 

Considering the homogeneous balance between $q''$ and $q^{3} $ in equation (\ref{9}), based on (\ref{12}) and (\ref{13}) we required that $m=1$. Taking into account the above considerations, the form of $q$ is as following

\begin{equation}
\label{14}
q(\xi )=a_{1} \left(\frac{G'}{G} \right)+a_{0}.
\end{equation} 

Then we know the exact view of $q^{3} $

\begin{equation}
\label{15}
q^{3} (\xi )=a_{1}^{3} \left(\frac{G'}{G} \right)^{3} +3a_{1}^{2} a_{0} \left(\frac{G'}{G} \right)^{2} +3a_{1} a_{0}^{2} \left(\frac{G'}{G} \right)+a_{0}^{3}.
\end{equation} 

Using (\ref{14}) and (\ref{11}) $q''$ is easily derived to

\begin{equation}
\label{16}
q''(\xi )=2a_{1} \left(\frac{G'}{G} \right)^{3} +3a_{1} \lambda \left(\frac{G'}{G} \right)^{2} +(2a_{1} \mu +a_{1} \lambda ^{2} )\left(\frac{G'}{G} \right)+a_{1} \lambda \mu .
\end{equation} 
By substituting (\ref{14})-(\ref{16}) into equation (\ref{9}) and collecting all terms with the same power of $(G'/G)$, the left-hand side of equation (\ref{9}) is converted into another polynomial in $(G'/G)$. 

\begin{equation}
\label{17}
\begin{aligned}
&(-2ka_{1}^{3} +2k^{3} a_{1} )\left(\frac{G'}{G} \right)^{3} +(-6ka_{1}^{2} a_{0} +3k^{3} a_{1} \lambda )\left(\frac{G'}{G} \right)^{2}\\&
+(\Omega '_{t} (t)a_{1} -6ka_{1} a_{0}^{2} +2k^{3} a_{1} \mu +k^{3} a_{1} \lambda ^{2} -k\gamma (t)q(0,t)a_{1} )\left(\frac{G'}{G} \right)\\&
+(C+\Omega '_{t} (t)a_{0} -2ka_{0}^{3} +a_{1} k^{3} \lambda \mu -k\gamma (t)q(0,t)a_{0} )=0.
\end{aligned}
\end{equation}
Equating each coefficient of expression (\ref{17}) to zero, yields a set of simultaneous equations for $a_{0} $, $a_{1} $, $\Omega (t)$ and $C$ which have the following form:

\begin{align*}
&\left(\frac{G'}{G} \right)^{3} : &\qquad &	-2ka_{1}^{3} +2k^{3} a_{1} =0,\\
&\left(\frac{G'}{G} \right)^{2} : &\qquad &  -6ka_{1}^{2} a_{0} +3k^{3} a_{1} \lambda =0,\\
&\left(\frac{G'}{G} \right)^{1} : &\qquad &  \Omega '_{t} (t)-6ka_{0}^{2} +2k^{3} \mu +k^{3} \lambda ^{2} -k\gamma (t)q(0,t)=0,\\
&\left(\frac{G'}{G} \right)^{0} : &\qquad &  C+\Omega '_{t} (t)a_{0} -2ka_{0}^{3} +a_{1} k^{3} \lambda \mu -k\gamma (t)q(0,t)a_{0} =0.
\end{align*}
By solving these the equations, we obtain the followings
\begin{equation}
\label{18}
\begin{aligned}
a_{0} &=\frac{k}{2} \lambda , \qquad  a_{1} =k, \qquad C=0,\\
\Omega (t)&=\frac{k^{3} (\lambda ^{2} -4\mu )}{2} t+k\int _{0}^{t}\gamma (\tau )q(0,\tau )d\tau +\Omega ^{0}.
\end{aligned}
\end{equation}
Here $\lambda $, $\mu $, $k$ and $\Omega ^{0} $ are arbitrary constants.
Using (\ref{18}), expression (\ref{14}) can be rewritten as
\begin{equation}
\label{19}
q(\xi )=k\left(\frac{G'}{G} \right)+\frac{k\lambda }{2},
\end{equation}
where $\xi =kx+\frac{k^{3}(\lambda ^{2} -4\mu )}{2} t+k\int _{0}^{t}\gamma (\tau )q(0,\tau )d\tau +\Omega ^{0}$. The function (\ref{19}) is a solution of equation (\ref{9}), provided that the integration constant $C$ in equation (\ref{9}) is taken as that in (\ref{18}). Substituting the general solutions of equation (\ref{11}) into (\ref{19}), we have three types of travelling wave solutions of the loaded mKdV equation (\ref{1}) as follows:

When $(\lambda ^{2} -4\mu )>0$,

\begin{equation}
\label{20}
q(\xi )=\frac{k\sqrt{\lambda ^{2} -4\mu } }{2} \left(\frac{c_{1} sh\frac{\sqrt{\lambda ^{2} -4\mu } }{2} \xi +c_{2} ch\frac{\sqrt{\lambda ^{2} -4\mu } }{2} \xi }{c_{1} ch\frac{\sqrt{\lambda ^{2} -4\mu } }{2} \xi +c_{2} sh\frac{\sqrt{\lambda ^{2} -4\mu } }{2} \xi } \right)
\end{equation}
where $\xi =kx+\frac{k^{3} (\lambda ^{2} -4\mu )}{2} t+k\int _{0}^{t}\gamma (\tau )q(0,\tau )d\tau +\Omega ^{0}$, $c_{1} $, $c_{2} $ and $\Omega ^{0} $ are arbitrary constants. It is obvious that the function $q(0,t)$ can be easily found based on expression (\ref{20}).

For example, let $\gamma(t)$ have a form below

$$\gamma (t)=\left(\frac{1}{k} \sum _{j=1}^{n}j\alpha _{j} t^{j-1} - \frac{k^{2} (\lambda ^{2} -4\mu )}{2} \right)\frac{2}{k\sqrt{\lambda ^{2} -4\mu } } cth\frac{\sqrt{\lambda ^{2} -4\mu } }{2} \sum _{j=0}^{n}\alpha _{j} t^{j},$$
where $\alpha _{j} \, (j=1,2,\, \, ...\, ,\, n)$ are constants. If  $c_{1} \ne 0$, $c_{2} =0$ and $(\lambda ^{2} -4\mu )>0$, then $q(x,t)$ becomes
\begin{equation}
\label{21}
q(x,t)=\frac{\sqrt{\lambda ^{2} -4\mu } }{2} k\, th\left(\frac{\sqrt{\lambda ^{2} -4\mu } }{2} (kx+\sum _{j=0}^{n}\alpha _{j} t^{j}  )\right).
\end{equation}
The function (\ref{21}) is the solution of the following loaded mKdV equation.

\begin{equation}
\label{21_1}
q_{t} -6q^{2} q_{x} +q_{xxx} -\left(\frac{1}{k} \sum _{j=1}^{n}j\alpha _{j} t^{j-1} - \frac{k^{2} (\lambda ^{2} -4\mu )}{2} \right)q_{x} =0.
\end{equation}

\begin{figure}[H]
	\centerline{%
		\includegraphics[width=18cm]{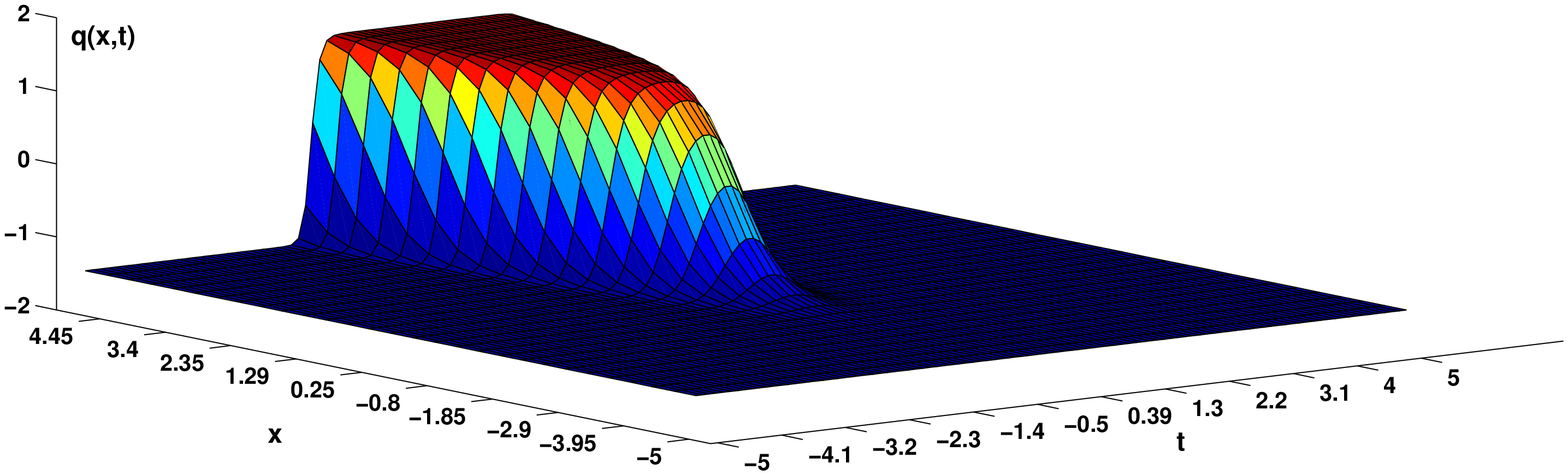}}
	\caption{The solution (\ref{21}) of the loaded mKdV equation (\ref{21_1}) for $\lambda =2$, $\mu =-1$, $k=1$, $\alpha_{0}=0$, $\alpha_{1}=1$, $\alpha_{2}=2$.}
\end{figure}
When $(\lambda ^{2} -4\mu )<0$, 
\begin{equation}
\label{22}
q(\xi )=\frac{k\sqrt{4\mu -\lambda ^{2} } }{2} \left(\frac{-c_{1} \sin \frac{\sqrt{4\mu -\lambda ^{2} } }{2} \xi +c_{2} \cos \frac{\sqrt{4\mu -\lambda ^{2} } }{2} \xi }{c_{1} \cos \frac{\sqrt{4\mu -\lambda ^{2} } }{2} \xi -c_{2} \sin \frac{\sqrt{4\mu -\lambda ^{2} } }{2} \xi } \right),
\end{equation}
where $\xi =kx-\frac{k^{3} (4\mu -\lambda ^{2} )}{2} t+k\int _{0}^{t}\gamma (\tau )q(0,\tau )d\tau +\Omega ^{0} $, $c_{1} $, $c_{2} $, and $\Omega ^{0} $ are arbitrary constants. It is not difficult for us to find $q(0,t)$ based on expression (\ref{22}). Let 

$$\gamma (t)=-\left(\frac{1}{k} \sum _{j=1}^{n}j\alpha _{j} t^{j-1} + \frac{k^{2} (4\mu -\lambda ^{2} )}{2} \right)\frac{2}{k\sqrt{4\mu -\lambda ^{2} } } ctg\left(\frac{\sqrt{4\mu -\lambda ^{2} } }{2} \sum _{j=0}^{n}\alpha _{j} t^{j}  \right),$$
where $\alpha _{j} \, (j=1,2,\, \, ...\, ,\, n)$ are constants, in particular, if  $c_{1} \ne 0$ and $c_{2} =0$, then $q(x,t)$ becomes

\begin{equation}
\label{23}
q(x,t)=-\frac{\sqrt{4\mu -\lambda ^{2} } }{2} k\, tg\left(\frac{\sqrt{4\mu -\lambda ^{2} } }{2} (kx+\sum _{j=0}^{n}\alpha _{j} t^{j} ) \right).
\end{equation}
The function (\ref{23}) is the solution of the following loaded mKdV equation.
\begin{equation}
\label{23_1}
q_{t} -6q^{2} q_{x} +q_{xxx} -\left(\frac{1}{k} \sum _{j=1}^{n}j\alpha _{j} t^{j-1} + \frac{k^{3} (4\mu -\lambda ^{2} )}{2} \right)q_{x} =0.
\end{equation}

\begin{figure}[H]
	\centerline{%
		\includegraphics[width=18cm]{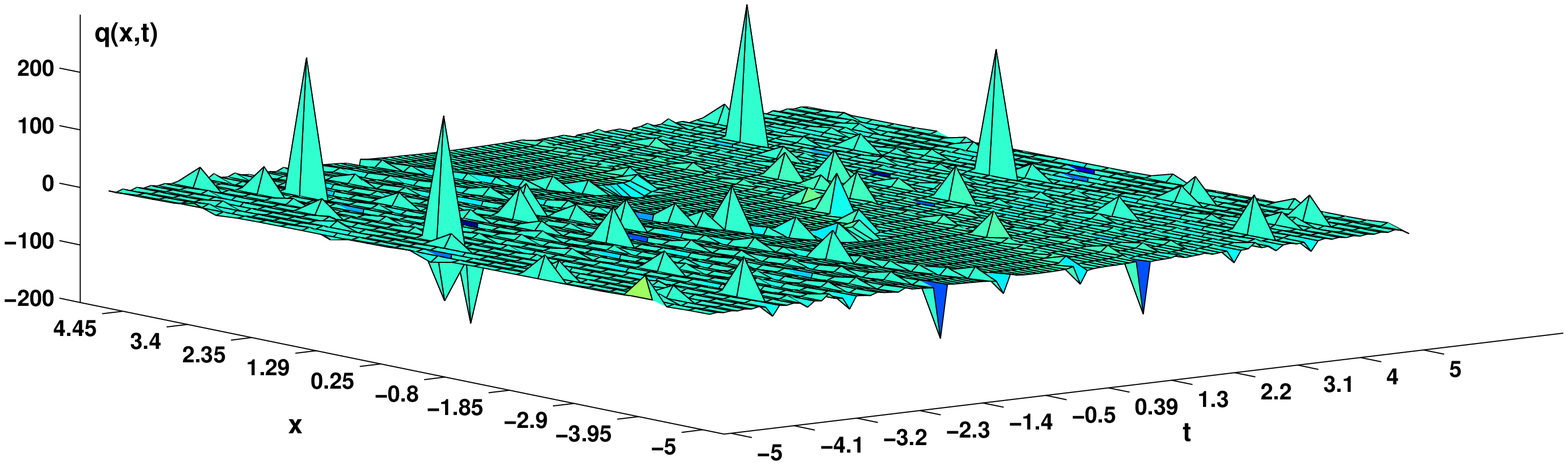}}
	\caption{The solution (\ref{23}) of the loaded mKdV equation (\ref{23_1}) for $\lambda =3$, $\mu =3$, $k=1$, $\alpha_{0}=0$, $\alpha_{1}=1$, $\alpha_{2}=2$.}
\end{figure}
When $(\lambda ^{2} -4\mu )=0$,

\begin{equation}
\label{24}
q(\xi )=\frac{kc_{2} }{c_{1} +\xi c_{2}},
\end{equation}
where $\xi =kx+k\int _{0}^{t}\gamma (\tau )q(0,\tau )d\tau +\Omega ^{0}  $, $c_{1} $, $c_{2} $, and $\Omega ^{0} $ are arbitrary constants. The function $q(0,t)$ is found based on expression (\ref{24}).
If  $c_{1} =0$, $c_{2} \ne 0$, $(\lambda ^{2} -4\mu )=0$ and $\gamma (t)$ has the form

$$\gamma (t)=\frac{\sum _{j=0}^{n}\alpha _{j} t^{j} \sum _{j=1}^{n}j\alpha _{j} t^{j-1}  }{k^{2} },$$
then $q(x,t)$ becomes
\begin{equation}
\label{25}
q(x,t)=\frac{k}{kx+\sum _{j=1}^{n}\alpha _{j} t^{j}}
\end{equation}
We know that the function (\ref{25}) satisfies the following loaded mKdV equation.
\begin{equation}
\label{25_1}
q_{t} -6q^{2} q_{x} +q_{xxx} -\left(\sum _{j=1}^{n}j\alpha _{j} t^{j-1}  \right)q_{x} =0.
\end{equation}

\begin{figure}[H]
	\centerline{%
		\includegraphics[width=18cm]{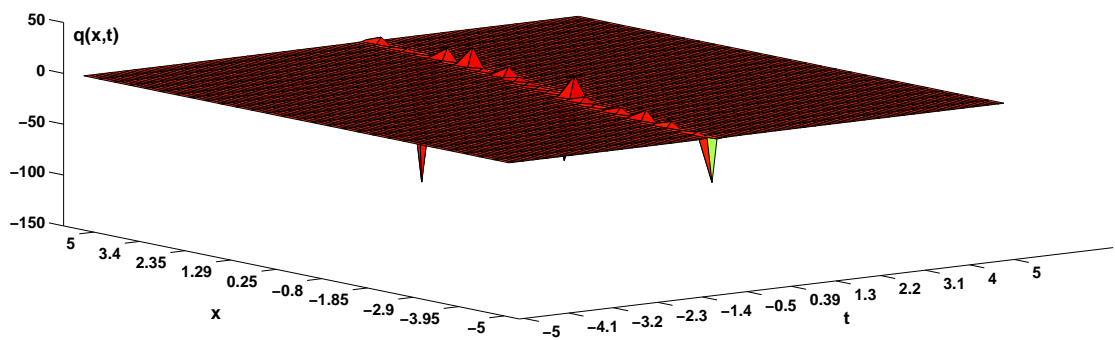}}
	\caption{The solution (\ref{25}) of the loaded mKdV equation (\ref{25_1}) for $\lambda =2$, $\mu =1$, $k=1$, $\alpha_{0}=0$, $\alpha_{1}=4$, $\alpha_{2}=1$, $\alpha_{3}=3$.}
\end{figure}

\section*{CONCLUSION}
The results of this study show that the $(G'/G)$ - extension method is effective in obtaining the exact solutions of the loaded modified Korteweg~--- de Vries equation. Parameters $c_ {1}$, $c_ {2}$, $\lambda$, $\mu$, $k$ and arbitrary function $\gamma(t)$ in solutions (\ref{20}), (\ref{20}), (\ref{24}) provide sufficient freedom for constructing solutions.

This method is easy to implement using well-known software packages, which allows you to solve complex nonlinear evolutionary equations of mathematical physics.

\end{document}